\documentclass[10pt]{article}

\usepackage[english]{babel}

\usepackage{amsfonts}
\usepackage{amssymb}

\newtheorem{theo}{Theorem}[section]
\newtheorem{prop}{Proposition}[section]
\newtheorem{lem}{Lemma}[section]
\newtheorem{cor}{Corollary}[section]
\newtheorem{DEF}{Definition}[section]
\newtheorem{EX}{Example}[section]
\newtheorem{REM}{Remark}[section]
\newenvironment{theorem}{\begin{theo}}{\end{theo}}  
\newenvironment{proposition}{\begin{prop}}{\end{prop}}  
\newenvironment{lemma}{\begin{lem}}{\end{lem}}
\newenvironment{corollary}{\begin{cor}}{\end{cor}}
\newenvironment{definition}{\begin{DEF}\rm}{\end{DEF}} 
\newenvironment{example}{\begin{EX}\rm}{\end{EX}} 
\newenvironment{remark}{\begin{REM}\rm}{\end{REM}}
\newenvironment{prof}{\noindent {\sl Proof.}}{\hfill $\square$ \par \ }

\newcommand{\gspazio}{\vskip1cm}
\newcommand{\mspazio}{\vskip.5cm}

\newcommand{\R}{\mathbb R}

\newcommand{\B}{\mathbb B}

\newcommand{\Hi}{\mathbb H}

\newcommand{\X}{\mathbb X}
\newcommand{\Y}{\mathbb Y}
\newcommand{\nullv}{\mathbf{0}}

\newcommand{\inte}{{\rm int}\, }

\newcommand{\fr}{{\rm bd}\, }
\newcommand{\dom}{{\rm dom}\, }
\newcommand{\graph}{{\rm gph}\, }

\newcommand{\Lin}{{\mathcal L}}

\newcommand{\der}{{\rm D}}
\newcommand{\Cone}{{\rm C}^{1}}
\newcommand{\Coneone}{{\rm C}^{1,1}}
\newcommand{\UCtwo}{{\rm UC}_2}
\newcommand{\CEff}{{\rm minimize}_C\ }

\newcommand{\modconv}[1]{\delta_{#1}}

\newcommand{\dist}[2]{{\rm dist}\left(#1, #2\right)}
\newcommand{\ball}[2]{{\rm B}\left(#1, #2\right)}
\newcommand{\reg}[2]{{\rm reg}\, #1({#2})}
\newcommand{\lip}[2]{{\rm lip}\,{#1}({#2})}
\newcommand{\Lip}[2]{{\rm Lip}({#1},{#2})}

\newcommand{\regat}[3]{{\rm reg}\, #1({#2}|{#3})}

\title{\bf Convexity of the images of small balls
through perturbed convex multifunctions}

\gspazio
\author{A. Uderzo \\
{\small\sl Dept. of Mathematics and Applications, University of Milano-Bicocca}\\
{\small\sl Via Cozzi, 55, 20125 Milano, Italy}\\
{\small e-mail:} {\small\tt amos.uderzo@unimib.it} }

\begin{document}

\maketitle

\vskip1cm

\begin{abstract} 
In the present paper, the following convexity principle is proved:
any closed convex multifunction, which is metrically regular in a
certain uniform sense near a given point, carries small balls
centered at that point to convex sets, even if it is perturbed by
adding $\Coneone$ smooth mappings with controlled Lipschizian
behaviour.
This result, which is valid for mappings defined on a subclass of
uniformly convex Banach spaces, can be regarded as a set-valued
generalization of the Polyak convexity principle. The latter,
indeed, can be derived as a special case of the former.
Such an extension of that principle enables one to build large
classes of nonconvex multifunctions preserving the convexity
of small balls. Some applications of this phenomenon to the theory
of set-valued optimization are proposed and discussed.
\end{abstract}

\mspazio

{\bf Keywords:} convex multifunction, uniformly convex Banach space,
modulus of convexity, metric regularity, Polyak convexity principle,
set-valued optimization.

\mspazio

{\bf Mathematics Subject Classification (AMS 2010):} 49J53, 52A05, 90C48.

\vfill\eject

%%%%%%%%%%%%%%%%%%%%%%%%%%%%%%%%%%%%%%%%%%%%%%%%%%%%%%%%%%
%%%%%%%%%%%%%%%%%%%%%%%%%%%%%%%%%%%%%%%%%%%%%%%%%%%%%%%%%%

\section{Introduction}

Treating problems from mathematical programming,
optimal control and from several areas of mathematical
economics yields a tremendous demand of convexity.
Convexity assumptions on problem data often strenghten
the analysis tools and trigger the application of special
approaches, otherwise not practicable. Even though
such a demand has led to deepen our knowledge
about convexity and then to develop expanding
branches of convex analysis, many fundamental
issues about convexity still remain to be investigated.
In the author's opinion, one of such issues concerns the
behaviour of convex sets under nonlinear transformations.
Indeed, not much seems to be known so far
about those sets whose image through nonlinear
mappings is convex. The existing results on this question can
be schematically classified as ``around a point" (local) results
or as ``on set" (nonlocal) results.
As an example of nonlocal result the Lyapunov
convexity theorem on the range of a vector measure
occupies a prominent place (see \cite{Lyap40}). It found
notable applications in control theory and mathematical
economics (see \cite{AliBor06,Olec90}). Other examples
of global results are, for instance, those in \cite{BoEmKo04,Reis07,Reis11}.
As an example of local result, the Polyak convexity principle
is certainly to be mentioned (see \cite{Poly01,Poly03}).
Like the Lyapunov's theorem, it revealed to be useful
in several topics of optimization and control theory, by
providing conditions upon which nonlinear mappings carry
small balls around a point to convex sets.

The present paper aims at bringing some contributions
in the same vein as the Polyak convexity principle, but entering
now the realm of set-valued mappings.
The starting point of the analysis here proposed is the
well-known fact that convex multifunctions (i.e. set-valued
mappings with convex graph) carry any convex set to
a convex set. If considering the category, whose objects
are convex sets, this class of mappings  seem to naturally
play the role of category morphisms. Unfortunately, by simple
examples it is readily realized that, when adding a nonlinear
single-valued mapping to a convex multifunction, in general
the convex graph property of the latter is broken.
Thus the question arises under which conditions mappings,
obtained by perturbing convex multifunctions by nonlinear
mappings, still carry small balls to convex sets. The main result
of this paper provides an answer to this problem.
It states that, if to a convex multifunction, which is metrically
regular near a reference point uniformly over its image,
a $\Coneone$ mapping is added, whose Lipschitzian behaviour
is controlled by the modulus of regularity of the former,
then the resulting set-valued mapping preserves
the convexity of small balls around the reference point of
its domain. In fact, this result can be regarded as an extension
of the Polyak convexity principle to a large class of
set-valued mappings. As it happens for its single-valued
counterpart, it is valid for mappings defined on uniformly
convex Banach spaces having second order polynomial
modulus of convexity. This class of spaces includes, for
instance, all Hilbert spaces. The proof combines a nice property,
coming from the rotund geometry of balls in the aforementioned
class of Banach spaces, with a convex solvability behaviour
of set-valued mappings, that are perturbed as described. The latter is a
consequence of the persistence of metric regularity under
additive Lipschitz perturbations, a well-known phenomenon
in variational analysis, which has revealed to be useful in various
contexts related to the solution stability and sensitivity for
generalized equations (see \cite{DonRoc09,Mord06}).

The contents of the paper are organized as follows. In Section
\ref{Sect:2} some tools, mainly from geometric functional analysis
and from nonlinear analysis, that are needed for establishing the
main result are recalled. In particular, in Subsection \ref{Sect:2.3}
a strenghtened notion of metric regularity for set-valued
mappings is introduced. Several classes of multifunctions
satisfying such a special property are exhibited, while it
is observed that the original notion of metric regularity
is weaker (in the sense that it holds more generally).
In Section \ref{Sect:3} the main result is proved and
commented. Then, it is shown how from this wider
convexity principle the Polyak's one can be derived,
as a special case.
Section \ref{Sect:4} is reserved to illustrate an application
of the main result to a topic from set-valued optimization.
More precisely, a class of optimization problems
is considered, whose set-valued objective is expressed
as a sum of a single-valued and a set-valued mapping.
This structure in the objective mapping may
model noise effects on vector optimization
problems. In this context, the convexity principle, under certain
additional assumptions, leads first of all
to establish the existence of efficient pairs for localizations
of an unconstrained problem, and then to achieve optimality
conditions based on the Lagrangian scalarization.

\section{Tools from nonlinear analysis}     \label{Sect:2}

\subsection{Uniformly convex Banach spaces}

The analysis of the posed problem will be carried out in the
particular setting of the uniformly convex real Banach spaces.
This because the main result presented in the paper essentially
rely on certain geometrical features of this specific class of
Banach spaces, features that are related to the rotundity of the
balls. The rotundity property of a ball in a Banach space $(\X,\|\cdot\|)$
can be quantitatively described by means of the function
$\modconv{\X}: [0,2]\longrightarrow [0,1]$, defined by
$$
    \modconv{\X}(\epsilon)=\inf\left\{1-\left\|\frac{x_1+x_2}{2}\right\|:\ 
    x_1,\, x_2\in\B,\ \|x_1-x_2\|\ge\epsilon\right\},
$$
which is called the {\em modulus of convexity} of $(\X,\|\cdot\|)$
{}\footnote{Equivalent definitions of the modulus of convexity
can be found in \cite{FaHaMoPeZi01}.}. 
$\B$ stands for the closed unit ball, centered at the null vector
$\nullv$ of $\X$.
Notice that $\modconv{\X}$ is not invariant under equivalent
renormings of $\X$. Such a notion allows one to define the class of
uniformly  convex Banach spaces, whose introduction is due to J.A.
Clarkson (see, for instance, \cite{Clar36,FaHaMoPeZi01,Megg98}). 

\begin{definition}
A Banach space  $(\X,\|\cdot\|)$ is called {\em uniformly convex}
(or, {\em uniformly rotund}) if
it is $\modconv{\X}(\epsilon)>0$ for every $\epsilon\in (0,2]$.
\end{definition}

In what follows, the modulus of convexity of a (uniformly convex)
Banach space is said to be of the {\em (polynomial) second order}
if there exists $c>0$ such that
$$
    \modconv{\X}(\epsilon)\ge c\epsilon^2,\quad\forall
   \epsilon\in [0,2].
$$
The class of uniformly convex real Banach spaces with second order
modulus of convexity reveals to be the proper setting, in which to
develp the analysis of the problem at the issue. Throughout the
paper, this class will be indicated by $\UCtwo$.

\begin{example}
($e_1$) By means of elementary considerations, the modulus of convexity
of a Hilbert space $\Hi$ can be calculated  to amount to
$$
    \modconv{\Hi}(\epsilon)=1-\sqrt{1-\frac{\epsilon^2}{4}},\quad\forall
    \epsilon\in [0,2].
$$
Therefore, every Hilbert space is uniformly convex, with a
second order modulus of convexity, such that $0<c\le 1/8$,
i.e. belongs to the class $\UCtwo$.

($e_2$) More generally, such Banach spaces as  $l^p$, $L^p$, and $W^p_m$,
with $1<p<2$, are known to have a modulus of convexity satisfying
the relation
$$
    \modconv{l^p}(\epsilon)=\modconv{L^p}(\epsilon)=\modconv{W^p_m}
    (\epsilon)> \frac{p-1}{8}\epsilon^2,\quad\forall\epsilon\in (0,2].
$$
Therefore, they also are examples of spaces of class $\UCtwo$
(see, for instance, \cite{FaHaMoPeZi01}).
\end{example}

\begin{remark}     \label{rem:UCspace}
($r_1$) Concerning the notion of uniform convexity, a caveat is due:
even finite-dimensional Banach spaces may fail to be uniformly
convex. Consider, for instance, $\R^2$ equipped with the Banach space
structure given by the norm $\|\cdot\|_\infty$.

($r_2$) It was proved that the modulus of convexity
$\modconv{\X}$ of any real Banach space, having dimension greater
than $1$, admits the following estimate from above
$$
    \modconv{\X}(\epsilon)\le1-\sqrt{1-\frac{\epsilon^2}{4}},\quad\forall
   \epsilon\in [0,2].
$$
This implies that the second order polynomial is a maximal one.

($r_3$) Recall that, according to the Milman-Pettis theorem,
every uniformly convex Banach space is reflexive, but the
converse is false (see, for instance, \cite{FaHaMoPeZi01}).
\end{remark}

For further material about uniformly convex Banach spaces,
see \cite{FaHaMoPeZi01, Megg98}.
In the following lemma, whose proof can be found in \cite{Uder13}
(Lemma 2.4), a key property of balls in any uniformly convex Banach
space of class $\UCtwo$ is stated, in view of a subsequent
application. Throughout the paper, given an element $x\in\X$ and
a real $r\ge 0$, $\ball{x}{r}$ denotes the closed ball centered at
the point $x$, with radius $r$.

\begin{lemma}     \label{lem:uniconv}
Let $(\X,\|\cdot\|)$ belong to $\UCtwo$, with modulus of convexity
$\modconv{\X}(\epsilon)\ge c\epsilon^2$, for some $c>0$.
Then, for every $x_0,\, x_1,\, x_2\in\X$ and
$r>0$, with $x_1,\, x_2\in\ball{x_0}{r}$, it holds
$$
    \ball{\frac{x_1+x_2}{2}}{\frac{c\|x_1-x_2\|^2}{r}}\subseteq
    \ball{x_0}{r}.
$$
\end{lemma}

%%%%%%%%%%%%%%%%%%%%%%

\vskip1cm

\subsection{Smooth mappings and Lipschitzian properties} 

Let $f:\Omega\longrightarrow\Y$ be a mapping between real Banach
spaces, where $\Omega$ is a nonempty open subset of $\X$. Its
G\^ateaux derivative at $\bar x\in\Omega$ is denoted by $\der f(\bar x)$.
Let us indicate by $(\Lin(\X,\Y),\|\cdot\|_\Lin)$ the Banach space of all
linear bounded operators between $\X$ and $\Y$, equipped with the
operator norm. If $f$ admits G\^ateaux derivative at each point of
 $\Omega$ and the mapping $\der f:\Omega\longrightarrow\Lin(\X,\Y)$,
defined by $x\mapsto \der f(x)$, is norm-to-$\|\cdot\|_\Lin$ continuous,
then $f$ is said to be of class ${\rm C}^1(\Omega)$. Remember that
if $f\in {\rm C}^1(\Omega)$, $f$ is in particular strictly differentiable
at each point of $\Omega$. If, furthermore, the mapping $\der f$ is Lipschitz
continuous on $\Omega$, $f$ is said to be of class $\Coneone
(\Omega)$. In such a case, the infimum of all constants $\kappa>0$
such that
$$
    \|\der f(x_1)-\der f(x_2)\|_\Lin\le\kappa\|x_1-x_2\|,\quad\forall
   x_1,\, x_2\in\Omega,
$$
will be indicated by $\Lip{\der f}{\Omega}$.
In the same setting, given a point $\bar x\in\Omega$, let us
define the value
$$
   \lip{f}{\bar x}=\limsup_{u,x\to\bar x\atop u\ne x}
   {\|f(u)-f(x)\|\over \|u-x\|}
$$
the Lipschitz modulus of $f$ at $\bar x$.
Clearly, $\lip{f}{\bar x}<\infty$ iff $f$ is locally Lipschitz in a
neighbourhood of $\bar x$. In particular, if $f\in {\rm C}^1(\inte
\ball{\bar x}{r})$ for some $r>0$, then one has $\lip{f}{\bar x}=\|\der
f(\bar x)\|_\Lin<\infty$. Throughout the paper, the convention is
adopted that, whenever $\|\der f(\bar x)\|_\Lin=0$ or $\lip{f}{\bar x}=0$,
the symbols $\|\der f(\bar x)\|_\Lin^{-1}$ and $\lip{f}{\bar x}^{-1}$
stand for $+\infty$.

This short subsection is concluded by a lemma, stating an estimate for
$\Coneone$ smooth mappings that will be crucially employed
in the proof of the main result. For its proof, the reader is
referred to \cite{Uder13} (Lemma  2.7).

\begin{lemma}     \label{lem:quadestim}
Let $f:\X\longrightarrow\Y$ be a mapping between Banach spaces,
let $U\subseteq\X$, let $\Omega\subseteq\X$ be an open set
such that $\Omega\supseteq U$, and let $x_1,\, x_2\in U$, with
$[x_1,x_2]\subseteq U$. If $f\in\Coneone(\Omega)$, then it holds
\begin{eqnarray*}
    \left\|\frac{f(x_1)+f(x_2)}{2}-f\left(\frac{x_1+x_2}{2}\right) \right\|\le
    \frac{\Lip{\der f}{U}}{8}\|x_1-x_2\|^2.
\end{eqnarray*}
\end{lemma}

%%%%%%%%%%%%%%%%%%%%%%

\vskip1cm

\subsection{Convex multifunctions and their metric
regularities}   \label{Sect:2.3}

Throughout the paper, given a subset $A$ of a Banach space
and a point $x$ in the same space, $\dist{x}{A}=\inf_{a\in A}\|a-x\|$
denotes the distance of $x$ from $A$.
The notion of metric regularity, along with its equivalent reformulations,
is recognized as an important tool in the variational analysis of
set-valued mappings. Recall that, given a set-valued mapping $G:
\X\rightrightarrows\Y$ between real Banach spaces, $G$ is said to
be metrically regular at $\bar x$, for $\bar y$, with $(\bar x,\bar y)
\in\graph G=\{(x,y)\in\X\times\Y:\ y\in G(x)\}$,
provided that there exist positive constants
$\kappa$, $\delta$, and $\zeta$ such that
\begin{eqnarray}      \label{in:mrdef}
      \dist{x}{G^{-1}(y)}\le\kappa\,\dist{y}{G(x)},\quad\forall
    x\in\ball{\bar x}{\delta},\forall y\in\ball{\bar y}{\zeta}.
\end{eqnarray}
The constant
$$
     \regat{G}{\bar x}{\bar y}=\inf\{\kappa\in (0,+\infty):\ (\ref{in:mrdef})
     \hbox{ holds for some $\delta$ and $\zeta$}\}
$$
is usually called regularity modulus of $G$ at $\bar x$, for $\bar y$.
Several aspects of the theory of metric regularity are exposed
in recent monographs (among the others, see \cite{BorZhu05,DonRoc09,
KlaKum02,Mord06,RocWet98}).

In what follows, a metric regularity property, which is stronger than the
original metric regularity at a reference pair, will be needed. Below,
given a real $r\ge 0$ and a subset $A\subset\Y$, by $\ball{A}{r}=\{y\in\Y:
\ \dist{y}{A}\le r\}$ the $r$-enlargement of $A$ will be indicated.

\begin{definition}       \label{def:mrset}
A set-valued mapping $G:\X\rightrightarrows\Y$ is said to be {\it metrically
regular at $\bar x\in\dom G=\{x\in\X:\ G(x)\ne\varnothing\}$,
for $G(\bar x)$}, if there exist positive
constants $\kappa$, $\delta$, and $\zeta$ such that
\begin{eqnarray}      \label{in:mrsetdef}
      \dist{x}{G^{-1}(v)}\le\kappa\,\dist{v}{G(x)},\quad\forall
    x\in\ball{\bar x}{\delta},\forall v\in\ball{G(\bar x)}{\zeta}.
\end{eqnarray}
The constant
$$
     \reg{G}{\bar x}=\sup_{y\in G(\bar x)}\regat{G}{\bar x}{y}
$$
will be used as a {\it regularity modulus} of $G$ at $\bar x$, for
$G(\bar x)$.
\end{definition}

From $(\ref{in:mrsetdef})$ one immediately sees that metric
regularity at $\bar x$, for $G(\bar x)$, implies (and is actually
equivalent to) the metric regularity of $G$ at $\bar x$, for every
$y\in G(\bar x)$, with the same constants $\kappa$, $\delta$,
and $\zeta$ in $(\ref{in:mrdef})$.
To the contrary, metric regularity at each pair
$\bar x$ and $y\in G(\bar x)$, without uniformity on the
values of $\kappa$, $\delta$, and $\zeta$ fails in general to
imply metric regularity of $G$ at $\bar x$, for $G(\bar x)$.
The example below illustrates such an occurence.

\begin{example}
Consider the function $g:\R^2\longrightarrow\R$, defined by
$$
    g(y_1,y_2)=y_1y_2,
$$
and, as a multifunction $G:\R\rightrightarrows\R^2$, its inverse
mapping
$$
    G(x)=g^{-1}(x)=\{y=(y_1,y_2)\in\R^2:\ y_1y_2=x\}.
$$
Set $\bar x=0$ and $\bar y=(\bar y_1,\bar y_2)=(0,0)$. Clearly,
$G(0)=\{y\in\R^2:\ y_1y_2=0\}$ is represented in the Euclidean plane
as the union of the two coordinate
axes. Observe that, since $g\in\Cone(\R^2)$, then $g$ is locally
Lipschitz near each point $y\in G(0)$. By consequence, according
to Theorem 1.49 in \cite{Mord06}, its inverse mapping $G$
turns out to be metrically regular at $0$, for each $y\in G(0)$.
Now, let $\kappa$, $\delta$ and $\zeta$ be arbitrary, but fixed,
positive reals. One has
$$
    \ball{G(0)}{\zeta}=[\R\times (-\zeta,\zeta)]\cup
    [(-\zeta,\zeta)\times\R].
$$
Notice that, if $x\in (-\delta,\delta)$ is close enough to $0$,
it is
$$
   G(x)\subseteq\ball{G(0)}{\zeta}.
$$
Let $x_\delta>0$ be such a point. Then, if $v=(v_1,v_2)\in
\ball{G(0)}{\zeta}\cap\R^2_+$, with $\R^2_+$ denoting the nonnegative
cone in $\R^2$, one sees that
$$
    \dist{v}{G(x_\delta)}<\zeta.
$$
Thus, choose $\bar v_2=\zeta/2$ and $\bar v_1$ in such a way that
$\bar v_1\bar v_2-x_\delta>\kappa\zeta$, i.e.
$$
    \bar v_1>{2(\kappa\zeta+x_\delta)\over \zeta}.
$$
It remains true that $\bar v\in\ball{G(0)}{\zeta}$, but one finds
$$
    \dist{x_\delta}{G^{-1}(\bar v)}=|x_\delta-\bar v_1\bar v_2|>
    \kappa\zeta>\kappa\dist{v}{G(x_\delta)}.
$$
So, inequality $(\ref{in:mrsetdef})$ is clearly violated.
\end{example}

Nonetheless, under additional assumptions on $G$, the metric regularity
at $\bar x$, for each point of $G(\bar x)$, can imply the metric regularity
at $\bar x$, for $G(\bar x)$. This happens, for instance, with
multifunctions taking compact values, as established in the next
proposition. Throughout the paper, given a subset $A$ of a Banach
space, by $\inte A$ the (topological) interior of $A$ is denoted.

\begin{proposition}     
Let $G:\X\rightrightarrows\Y$ be a set-valued mapping between
Banach spaces, and let $\bar x\in\dom G$. If $G$ is metrically
regular at $\bar x$, for every $y\in G(\bar x)$, and $G(\bar x)$
is compact, then $G$ is metrically regular at $\bar x$, for
$G(\bar x)$.
\end{proposition}

\begin{prof}
By virtue of the metric regularity of $G$ at $\bar x$, for every
$y\in G(\bar x)$, there exist positive $\delta_y$, $\zeta_y$ and
$\kappa_y$ such that
\begin{eqnarray}       \label{in:mratbarxy}
    \dist{x}{G^{-1}(v)}\le\kappa_y\, \dist{v}{G(x)},\quad\forall
   x\in\ball{\bar x}{\delta_y},\ \forall v\in\ball{y}{\zeta_y}.
\end{eqnarray}
Notice that the family $\{\inte\ball{y}{\zeta_y/2}:\ y\in G(\bar x)\}$
forms an open covering of $G(\bar x)$. Since $G(\bar x)$ has
been supposed to be compact, this family must admit a finite
subfamily still covering $G(\bar x)$, say $\{\inte\ball{y_i}{\zeta_{y_i}/2}:
\ y_i\in G(\bar x),\ i=1,\dots, m\}$. Thus, it is possible to define
the following positive constants
$$
   \delta=\min\{\delta_{y_i}:\ i=1,\dots,m\},\qquad
   \zeta=\min\{\zeta_{y_i}:\ i=1,\dots,m\},\qquad\hbox{and}
   \qquad\kappa=\max\{\kappa_{y_i}:\ i=1,\dots,m\}.
$$
Now, if $v\in\ball{G(\bar x)}{\zeta/3}$, there must exist $y\in
G(\bar x)$ such that $d(v,y)<\zeta/2$. Since it is
$$
    y\in G(\bar x)\subseteq\bigcup_{i=1}^m\inte\ball{y_i}{\zeta_{y_i}/2},
$$
then for some index $i^*\in\{1,\dots,m\}$ one has $y\in\inte
\ball{y_{i^*}}{\zeta_{y_{i^*}}/2}$. It follows
$$
   d(v,y_{i^*})\le d(v,y)+d(y,y_{i^*})<{\zeta\over 2}+{\zeta_{y_{i^*}}\over 2}
   \le\zeta_{y_{i^*}}.
$$
Hence it is possible to invoke inequality $(\ref{in:mratbarxy})$, in the case
$y=y_{i^*}$. Consequently, one obtains
$$
     \dist{x}{G^{-1}(v)}\le\kappa\, \dist{v}{G(x)},\quad\forall
   x\in\ball{\bar x}{\delta}, \ v\in\ball{G(\bar x)}{\zeta/3}.
$$
This completes the proof.
\end{prof}

Further examples of multifunctions satisfying Definition \ref{def:mrset}
can be found within the class of convex multifunctions, that plays
a leading role in the present work. Let us recall that a set-valued
mapping $G:\X\rightrightarrows\Y$ between Banach spaces is said
to be convex if $\graph G$ is a convex set. Equivalently, $G$ is
convex iff
$$
    tG(x_1)+(1-t)G(x_2)\subseteq G(tx_1+(1-t)x_2),\quad\forall
    t\in [0,1],\quad\forall x_1,\, x_2\in\X,
$$
with the convention that $\varnothing+S=\varnothing=t\varnothing$,
for every $S\subseteq\Y$ and $t\in\R$ (see \cite{AubFra90}).

Whenever a convex multifunction $G:\X\rightrightarrows\Y$ is also
positively homogeneous, i.e.
$$
    \nullv\in G(\nullv) \qquad\hbox{ and }\qquad
   G(\lambda x)=\lambda G(x),\quad\forall\lambda>0,
   \forall x\in\X,
$$
it  is called sublinear or, according to \cite{AubFra90,Rock70,Robi72},
a convex process. In other terms, sublinear set-valued mappings
are characterized by having a cone in $\X\times\Y$ as their graph.
During the 70-ies and the 80-ies, they have been the subject
of deep investigations in
convex and nonsmooth analysis. In particular, the study of their
regularity properties has revealed that the value $\regat{G}
{\nullv}{\nullv}$ plays a crucial role in understanding their special
behaviour. More precisely, it is known that
$$
     \regat{G}{\nullv}{\nullv}=\|G^{-1}\|^-,
$$
where 
$$
   \|H\|^-=\sup_{x\in\B}\inf_{y\in H(x)}\|y\|=\sup_{x\in\B}
   \dist{\nullv}{H(x)},
$$
is the so-called inner
norm of a positively homogeneous set-valued mapping $H:\X
\rightrightarrows\Y$ (see \cite{DoLeRo03,DonRoc09,Mord06}).
Furthermore, it has been shown that for any sublinear mapping
$G:\X\rightrightarrows\Y$ with closed graph it results in
\begin{eqnarray}
    \regat{G}{\bar x}{\bar y}\le \regat{G}{\nullv}{\nullv},
   \quad\forall (\bar x,\bar y)\in\graph G
\end{eqnarray}
(see, for instance, \cite{DonRoc09}).
In terms of the metric regularity notion introduced in Definition
\ref{def:mrset}, such a property can be restated as follows.

\begin{proposition}    \label{pro:sublmr}
Let $G:\X\rightrightarrows\Y$ be a sublinear set-valued mapping
between real Banach spaces and let $(\bar x,\bar y)\in\graph G$.
If $\regat{G}{\nullv}{\nullv}<\infty$, then $G$ is metrically regular
at $\bar x$, for $G(\bar x)$, and it holds
\begin{eqnarray*}
    \reg{G}{\bar x}\le\regat{G}{\nullv}{\nullv}=\reg{G}{\nullv}.
\end{eqnarray*}
\end{proposition}

After the works of Lyusternik, Graves, Robinson and Milyutin,
it was well
understood that the regularity property is stable under additive
perturbations with locally Lipschitz mappings, provided that
the Lipschitz modulus is small enough. The following result
provides a quantitative description of such a persistence
phenomenon (see \cite{DoLeRo03,DonRoc09,Mord06}).

\begin{theorem}       \label{thm:DoLeRo}
{\bf (estimate for Lipschitz perturbations)}
Consider a mapping $G:\X\rightrightarrows\Y$ and $(\bar x,\bar y)\in
\graph G$, at which $\graph G$ is locally closed, and a mapping
$f:\X\longrightarrow\Y$. If $\regat{G}{\bar x}{\bar y}<\kappa<\infty$
and $\lip{f}{\bar x}<\lambda<\kappa^{-1}$, then
$$
    \regat{(f+G)}{\bar x}{f(\bar x)+\bar y}<{1\over\kappa^{-1}-\lambda}.
$$
\end{theorem}

In the next lemma, a uniform behaviour of the metric regularity
property as given in Definition \ref{def:mrset}
in the presence of additive Lipschitz perturbations is
obtained. It will be exploited in the proof of the main result.

\begin{lemma}      \label{lem:unimrpert}
Let $G:\X\rightrightarrows\Y$ be a set-valued mapping, with $\graph
G$ locally closed, and
let $f:\X\longrightarrow\Y$ and let $\bar x\in\dom G$.
Suppose that $G$ is metrically regular at $\bar x$, for $G(\bar x)$,
$f$ is locally Lipschitz near $\bar x$, and
\begin{eqnarray}      \label{in:lemunimrmod}
   \reg{G}{\bar x}<\lip{f}{\bar x}^{-1}.
\end{eqnarray}
Then, the set-valued mapping $F=f+G$ is metrically regular
at $\bar x$, for $F(\bar x)=f(\bar x)+G(\bar x)$. Moreover
\begin{eqnarray*}
   \reg{F}{\bar x}\le {1\over \reg{G}{\bar x}^{-1}
   -\lip{f}{\bar x}},\quad\forall y\in\ G(\bar x).
\end{eqnarray*}
\end{lemma}

\begin{prof}
Fix an arbitrary $y\in G(\bar x)$. By inequality $(\ref{in:mrsetdef})$,
taking an arbitrary $\kappa$, with $\kappa> \reg{G}{\bar x}\ge
\regat{G}{\bar x}{y}$, one obtains
\begin{eqnarray*}  
      \dist{x}{G^{-1}(v)}\le\kappa\,\dist{v}{G(x)},\quad\forall
    x\in\ball{\bar x}{\delta},\forall v\in\ball{y}{\zeta}.
\end{eqnarray*}
By proceeding as in the proof of Theorem \ref{thm:DoLeRo}
(Theorem 3.3 in \cite{DoLeRo03}), it is possible to find
values of $\tilde\kappa$, $\tilde\delta$ and $\tilde\zeta$,
depending only on $\kappa$, $\delta$ and $\zeta$ (but not on $y$!)
\footnote{In the proof of Theorem 3.3, the new constants
for which inequlity $(\ref{in:mrtildeF})$ holds are expressed
in terms of  $\kappa$, $\delta$ and $\zeta$ only.},
such that
\begin{eqnarray}  \label{in:mrtildeF}
      \dist{x}{F^{-1}(v)}\le\tilde\kappa\,\dist{v}{F(x)},\quad\forall
     x\in\ball{\bar x}{\tilde\delta},\forall v\in\ball{f(\bar x)+y}
    {\tilde\zeta},
\end{eqnarray}
with $\tilde\kappa=(\kappa-\lambda)^{-1}$, for any $\lambda\in
(\lip{f}{\bar x},1/ \reg{G}{\bar x})$.
The infimum over all values of $\tilde\kappa$ such that inequality
$(\ref{in:mrtildeF})$ holds true can be shown consequently
not to exceed $(\regat{G}{\bar x}{y}^{-1}-\lip{f}{\bar x})^{-1}$, and
hence the value $(\reg{G}{\bar x}^{-1}-\lip{f}{\bar x})^{-1}$.
Since if $v\in\ball{F(\bar x)}{\tilde\zeta}$, then a $y\in G(\bar x)$
must exist such that $v\in\ball{f(\bar x)+y}{\tilde\zeta}$, one
gets the validity of inequality $(\ref{in:mrsetdef})$.
According to the definition of $\reg{F}{\bar x}$, this completes
the proof.
\end{prof}

\vskip1cm

%%%%%%%%%%%%%%%%%%%%%%%%%%%%%%%%%%%%%%%%%%%%%%%%%%%%%%%%%%%

\section{The main result}     \label{Sect:3}

One is now in a position to establish the following sufficient
condition for the convexity of the images of small balls
through a convex multifunction $G$ perturbed by a $\Coneone$
mapping $f$, which is the main result of the paper.

\begin{theorem}     \label{thm:PCPset}
Let $G:\X\rightrightarrows\Y$ be a set-valued mapping between real Banach
spaces, let $f:\Omega\longrightarrow\Y$ be a mapping defined on an open
set $\Omega$ and let $x_0$ and
$r>0$ such that $\ball{x_0}{r}\subseteq\Omega\cap\dom G$. Suppose that:

\noindent (i) $(\X,\|\cdot\|)$ is of class $\UCtwo$, having second order
modulus of convexity with some constant $c>0$;

\noindent (ii) $f\in\Coneone(\inte\ball{x_0}{r})$;

\noindent (iii) $G$ is a closed and convex multifunction;

\noindent (iv) $G$ is upper semicontinuous (for short, u.s.c.)
at $x_0$;

\noindent (v) $G$ is metrically regular at $x_0$, for $G(x_0)$,
with regularity modulus such that
\begin{eqnarray}       \label{in:mrGlipf}
    \reg{G}{x_0}< \|\der f(x_0)\|_\Lin^{-1};
\end{eqnarray}
   
\noindent (vi)  there exists $\tau>0$ such that $F(\ball{x_0}{t})$
is closed for every $t\in [0,\tau]$.

\noindent Then, there exists $\epsilon_0>0$, such that
$F(\ball{x_0}{\epsilon})$ is convex, for every $\epsilon\in
[0,\epsilon_0]$.
\end{theorem}

\begin{prof}
As already remarked, since $f\in\Coneone(\inte\ball{x_0}{r})$,
then it is $\lip{f}{x_0}= \|\der f(x_0)\|_\Lin<\infty$. According to
hypothesis $(v)$, $G$ is metrically regular at $x_0$, for $G(x_0)$,
and condition $(\ref{in:lemunimrmod})$ takes place.
Thus, by virtue of Lemma \ref{lem:unimrpert}, the set-valued
mapping $F=f+G$ is metrically regular at $x_0$, for each $f(x_0)+y$,
with $y\in G(x_0)$, that is there exist $\delta>0$ and $\zeta>0$ such that
\begin{eqnarray}      \label{in:mrx0}
   \dist{x}{F^{-1}(v)}\le\kappa\,\dist{v}{F(x)},\quad\forall
   x\in\ball{x_0}{\delta},\ \forall v\in\ball{f(x_0)+y}{\zeta},
\end{eqnarray}
for any $\kappa>(\reg{G}{x_0}^{-1}- \|\der f(x_0)\|_\Lin)^{-1}$,
and it holds
$$
    \regat{F}{x_0}{f(x_0)+y}\le{1\over \reg{G}{x_0}^{-1}-
    \|\der f(x_0)\|_\Lin}.
$$
Recall that the constants appearing in inequality $(\ref{in:mrx0})$ remain the
same for every $y\in G(x_0)$.  Then, corresponding to $\zeta/4$,
as a consequence of hypothesis $(ii)$, by continuity of $f$
at $x_0$,  there is $\delta_1>0$
such that
$$
   f(x)\in\ball{f(x_0)}{\zeta/4},\quad\forall x\in\ball{x_0}{\delta_1}.
$$
Again, by upper semicontinuity of $G$ at $x_0$ (hypothesis
$(iv)$), corresponding to $\zeta/4$, there is $\delta_2>0$ such that
$$
   G(x)\subseteq\ball{G(x_0)}{\zeta/4},\quad\forall x\in
   \ball{x_0}{\delta_2}.
$$ 
Consequently, take $\epsilon_0$ in such a way that
\begin{eqnarray}
   0<\epsilon_0<\left\{\delta,\,\delta_1,\,\delta_2,\,\tau,\, r,\,
  {4c( \reg{G}{x_0}^{-1}- \|\der f(x_0)\|_\Lin)\over
   \Lip{\der f}{\inte\ball{x_0}{r}}+1}\right\}.
\end{eqnarray}
In the case $\epsilon=0$ the thesis becomes trivial, because
$F(x_0)=f(x_0)+G(x_0)$ is convex as a sum of the convex sets
$\{f(x_0)\}$ and $G(x_0)$ (the latter is convex as a consequence of
hypothesis $(iii)$).
Now, fix $\epsilon\in (0,\epsilon_0]$. Since it is $\epsilon<\tau$ and hence,
according to hypothesis $(vi)$ set $F( \ball{x_0}{\epsilon})$ is closed,
to show that this set is convex it suffices to prove that, whenever $y_1$,
$y_2\in F( \ball{x_0}{\epsilon})$, it happens also that
$$
    \bar y={y_1+y_2\over 2}\in F( \ball{x_0}{\epsilon}).
$$
The fact that $y_1\in F(\ball{x_0}{\epsilon})$ implies the existence
of $x_1\in\ball{x_0}{\epsilon}$ such that $y_1\in F(x_1)=f(x_1)+G(x_1)$,
and hence the existence of $v_1\in G(x_1)$ such that $y_1=f(x_1)+v_1$.
Analogously, the fact that $y_2\in F(\ball{x_0}{\epsilon})$ implies
the existence of $x_2\in\ball{x_0}{\epsilon}$ and $v_2\in G(x_2)$,
such that $y_2=f(x_2)+v_2$.
Set
$$
   \bar v={v_1+v_2\over 2}\qquad\hbox{ and }\qquad
    \bar x={x_1+x_2\over 2}.
$$
If $\bar y\in F(\bar x)\subseteq F(\ball{x_0}{\epsilon})$ the argument
is finished. Otherewise, it is $\dist{\bar y}{F(\bar x)}>0$ because
$F(\bar x)$ is closed.
Notice that, since it is $\epsilon<\delta_1$, one has $f(x_1)$, $f(x_2)\in
\ball{f(x_0)}{\zeta/4}$ and hence
\begin{eqnarray}     \label{inc:midvalf}
    {f(x_1)+f(x_2)\over 2}\in\ball{f(x_0)}{\zeta/4}.
\end{eqnarray}
Since it is  $\epsilon<\delta_2$, one has that  $v_1$, $v_2\in
\ball{G(x_0)}{\zeta/4}$. The fact that $G$ is a convex multifunction
implies that $G(x_0)$ is convex, and so is function $v\mapsto
\dist{v}{G(x_0)}$, with the consequence that
$$
   \bar v\in\ball{G(x_0)}{\zeta/4}.
$$
This means that there exists $y_0\in G(x_0)$ such that
$d(\bar v,y_0)<\zeta/2$.
Then, from inequality $(\ref{inc:midvalf})$ it follows
\begin{eqnarray*}
    d(\bar y,f(x_0)+y_0) &=& \left\|{f(x_1)+f(x_2)\over 2}+\bar v
    -(f(x_0)+y_0)\right\|\le
    \left\|{f(x_1)+f(x_2)\over 2}-f(x_0)\right\|+\|\bar v-y_0\|   \\
    &\le& {\zeta\over 4}+{\zeta\over 2}<\zeta.
\end{eqnarray*}
The above inequalities show that $\bar x\in\ball{x_0}{\delta}$
and $\bar y\in\ball{f(x_0)+y_0}{\zeta}$, so inequality $(\ref{in:mrx0})$
applies, namely for any $\kappa>(\reg{G}{x_0}^{-1}- \|\der
f(x_0)\|_\Lin)^{-1}$ it holds
$$
    \dist{\bar x}{F^{-1}(\bar y)}\le\kappa\, \dist{\bar y}{F(\bar x)}.
$$
As a consequence of the last inequality, there exists $\hat x
\in F^{-1}(\bar y)$ such that
\begin{eqnarray}     \label{in:hatbardist}
   d(\bar x,\hat x)<{2\,\dist{\bar y}{F(\bar x)}\over\reg{G}{x_0}^{-1}
   - \|\der f(x_0)\|_\Lin}.
\end{eqnarray}
Now, observe that, by an obvious translation of vectors, one
obtains
$$
   \dist{\bar y}{F(\bar x)}=\dist{{f(x_1)+f(x_2)\over 2}+\bar v}
   {f(\bar x)+G(\bar x)}=
   \dist{{f(x_1)+f(x_2)\over 2}-f(\bar x)}{G(\bar x)-\bar v}.
$$
Since, by convexity of $\graph G$, it is
$$
     \bar v\in {G(x_1)+G(x_2)\over 2}\subseteq G(\bar x),
$$
in the light of Lemma \ref{lem:quadestim} it results in
\begin{eqnarray*}
   \dist{\bar y}{F(\bar x)} &\le& \dist{{f(x_1)+f(x_2)\over 2}-f(\bar x)}
   {{G(x_1)+G(x_2)\over 2}-\bar v}\le
   \left\|{f(x_1)+f(x_2)\over 2}-f(\bar x)\right\| \\
   &\le&  {\Lip{\der f}{\inte\ball{x_0}{r}}\over 8}\|x_1-x_2\|^2.
\end{eqnarray*}
From inequality $(\ref{in:hatbardist})$, recalling that
$$
   \epsilon<  {4c( \reg{G}{x_0}^{-1}- \|\der f(x_0)\|_\Lin)\over
   \Lip{\der f}{\inte\ball{x_0}{r}}+1},
$$
one obtains
$$
   d(\hat x,\bar x)<{\Lip{\der f}{\ball{x_0}{r}}\over
  4( \reg{G}{x_0}^{-1}- \|\der f(x_0)\|_\Lin)}\|x_1-x_2\|^2<
  {c\over\epsilon}\|x_1-x_2\|^2,
$$
whence it follows that
$$
    \hat x\in\ball{\bar x}{{c\over\epsilon}\|x_1-x_2\|^2}.
$$
By the uniform convexity of $\X$, with modulus of second order
of constant $c$, the last inclusion is known to imply that
$\hat x\in\ball{x_0}{\epsilon}$, according to Lemma \ref{lem:uniconv}.
Thus
$$
   \bar y\in f(\hat x)+G(\hat x)\subseteq F(\ball{x_0}{\epsilon}).
$$
The arbitrariness of $\epsilon\in (0,\epsilon_0]$ completes
the proof.
\end{prof}

\begin{remark}     \label{rem:mainthm}
$(r_1)$ The reader should notice that Theorem \ref{thm:PCPset}
has a local nature. Therefore hypothesis $(iii)$, the only global
one, can be actually weakened by assuming $G$ to be locally closed
convex near $x_0$ and $G(x_0)$, i.e. that there exists $r>0$
such that $\graph G\cap [\ball{x_0}{r}\times \ball{G(x_0)}{r}]$
is closed and convex. A perusal of the arguments in the proof
confirms the validity of such a refinement.

$(r_2)$ Whenever $G:\X\rightrightarrows\Y$ is, in particular, a closed
sublinear set-valued mapping, then in the light of the global metric
regularity property recalled in Proposition \ref{pro:sublmr}, hypothesis
$(v)$ takes the simpler form: $G$ is metrically regular at $\nullv$,
for $\nullv$, and $\regat{G}{\nullv}{\nullv}<\|\der f(x_0)\|_\Lin^{-1}$.

$(r_3)$ As a consequence of the metric regularity of $F=f+G$ at
$x_0$, for $F(x_0)$, it follows that if $x\in\inte\ball{x_0}{\epsilon}$,
with $\epsilon\in (0,\epsilon_0]$, and $y\in F(x)$, then $y\in\inte
F(\ball{x_0}{\epsilon})$. Thus, if denoting by $\fr A$ the boundary
of a subset $A$, whenever $y\in\fr F(\ball{x_0}{\epsilon})$ and
$x\in F^{-1}(y)$, one obtains that $x\not\in\inte\ball{x_0}{\epsilon}$,
namely $x\in\fr\ball{x_0}{\epsilon}$. In particular, one
has that $\inte F(\ball{x_0}{\epsilon})\ne\varnothing$.

$(r_4)$ The strict inequality appearing in $(\ref{in:mrGlipf})$
is essential and can not be relaxed by a non strict one, even
in very simple cases, as illustrated by the counterexample below.
\end{remark}

\begin{example}
Let $f:\R\longrightarrow\R^2$ and $G:\R\rightrightarrows\R^2$ be
given by
$$
    f(x)=(0,\, x^2) \qquad\hbox{ and }\qquad G(x)=\{(x,x)\},
$$
respectively, and let $x_0=0$, with $\R$ and $\R^2$ equipped with
their usual Euclidean structure. Then, it results in
$$
   F(x)=f(x)+G(x)=\{(x,x^2+x)\}.
$$
Notice that $f\in\Coneone(\R)$ and $G$ is a convex process.
Throught elementary calculations, one finds $\|\der f(0)\|=0$
and, since $G$ is not onto, $\reg{G}{0}=+\infty$. In other words,
the stric inequality $(\ref{in:mrGlipf})$ is not true, being
replaced by an equality. As one
easily checks, all remaining hypotheses of Theorem \ref{thm:PCPset}
are fulfilled. In this case the thesis fails to be true. Indeed, the
image of a ball $\ball{0}{\epsilon}=[-\epsilon,\epsilon]$ through
$F$ is the set
$$
    F([-\epsilon,\epsilon])=\{(x,x^2+x)\in\R^2:\ -\epsilon\le x\le\epsilon\},
$$
that fails to be convex, for every $\epsilon>0$.
\end{example}

From Theorem \ref{thm:PCPset} one can derive, as a special case,
a sufficient condition for the convexity of images of small balls,
around a regular point, which is known as a Polyak's convexity
principle.

\begin{corollary} {\bf (Polyak convexity principle)}  \label{thm:Polyakteo}
Let $f:\X\longrightarrow\Y$ be a mapping between real Banach
spaces, let $\Omega$ be an open subset of $\X$, let $x_0\in\Omega$,
and $r>0$ such that $\ball{x_0}{r}\subseteq\Omega$. Suppose that:

\noindent $(i)$ $(\X,\|\cdot\|)$ is of class $\UCtwo$;

\noindent $(ii)$ $f\in\Coneone(\Omega)$ and $\der f(x_0)\in
\Lin(\X,\Y)$ is onto.

\noindent Then, there exists $\epsilon_0\in (0,r)$ such that
$f(\ball{x_0}{\epsilon})$ is convex, for every $\epsilon\in
[0,\epsilon_0]$.
\end{corollary}

\begin{prof}
Observe that, under the current hypotheses, the mapping
$x\mapsto\{\der f(x_0)[x]\}$ is a closed sublinear multifunction,
which is u.s.c. at $x_0$, as $\der f(x_0)\in\Lin(\X,\Y)$. According
to the Banach-Schauder theorem, the fact that $\der f(x_0)$ is
onto is equivalent to its global metric regularity, and it holds
$$
    \regat{\der f(x_0)}{\nullv}{\nullv}=\|\der f(x_0)^{-1}\|^-<\infty
$$
(here $\der f(x_0)^{-1}$ denotes the multivalued inverse of
$\der f(x_0)$). Therefore, it remains to set
$$
    h=f-\der f(x_0),
$$
so that $f=h+\der f(x_0)$ can be expressed as a perturbation
of $\der f(x_0)$. Clearly $h\in\Coneone(\inte\ball{x_0}{r})$ and
$\der h(x_0)=\nullv\in\Lin(\X,\Y)$, so, according to the convention
made, condition $(\ref{in:mrGlipf})$ is fulfilled, independently of the
value of $\|\der f(x_0)^{-1}\|^-$. Finally, in Lemma 2.10 of \cite{Uder13}
the closedness of $f(\ball{x_0}{t})$, for every $t\in [0,\tau]$, has
been shown to come as a consequence of the metric regularity
of $f$ at $x_0$, which is in turn a consequence of the surjectivity
of $\der f(x_0)$, as it is know by the Lyusternik-Graves theorem.
Thus, Theorem \ref{thm:PCPset} applies.
\end{prof}

%%%%%%%%%%%%%%%%%%%%%%%%%%%%%%%%%%%%%%%%%%%%%%%%%%%%%%%%%

\section{An application to set-valued optimization}     \label{Sect:4}

In this section an application of the main result is presented,
which concerns set-valued optimization. This is a rather recent
branch of optimization, focusing on problems whose objective
(or cost) function are set-valued mappings. Some motivating
examples, coming from applications to mathematical economics
as well as from theoretical issues in vector optimization, fuzzy
programming and robust optimization, are described, for instance,
in \cite{BaoMor10,BaoMor11,KhTaZa15}.

In what follows, let us assume that a vector objective function
$q:\X\longrightarrow\Y$, acting in abstract spaces, is given
as a problem datum, along with a partial ordering $\le_C$
on its range space, which is defined by a proper, convex,
pointed and closed cone $C\subset\Y$. In real-world scenarios,
it may happen that the value of $q$ is affected by noise
effects, due to approximations, errors and/or incompleteness
in measurement and informations. As a result, instead of a unique
vector cost $q(x)$ corresponding to a chosen strategy $x$ in the decision
space $\X$, one has to deal with a set of several vectors in $\Y$.
This situation can be formalized by assuming that $q$ is perturbed
by adding a given set-valued mapping $Q:\X\rightrightarrows\Y$,
leading to a set-valued objective $\Phi=q+Q$. The resulting
(unconstrained) optimization problem is
$$
   \CEff \Phi(x) \quad \hbox{ over }\Omega,
    \leqno (\mathcal{SP})
$$
where $\Omega$ is a nonempty open subset of $\X$. Throughout
the present section, it will be assumed that $\dom\Phi\supseteq
\Omega$.

For such a problem several solution concepts have been proposed.
Following a vector based approach, according to \cite{KhTaZa15}
a pair $(\bar x,\bar y)\in\graph\Phi$ is said to be a $C$-efficient
pair for problem $(\mathcal{SP})$ if
$$
    (\bar y-C)\cap\Phi(\Omega)=\{\bar y\}.
$$
Notice that $\bar y$ is a $C$-minimal element of $\Phi(\Omega)$
with respect to the partial order relation $\le_C$.
In this context, as a consequence of Theorem \ref{thm:PCPset},
 the existence of $C$-efficient pairs of localizations of problem
$(\mathcal{SP})$ is established. Given a point $x_0\in\Omega$
and $\epsilon>0$, by a localization of problem $(\mathcal{SP})$
the following constrained set-valued minimization
problem is meant
$$
   \CEff \Phi(x)\quad \hbox{ subject to } x\in\ball{x_0}{\epsilon}.
   \leqno (\mathcal{SP}_{x_0,\epsilon})
$$

The following technical lemma will be employed in the proof of the
next result.

\begin{lemma}     \label{lem:cuscblocb}
Let $q:\X\longrightarrow\Y$ and $Q:\X\rightrightarrows\Y$
be given. Suppose that $q$ is continuous at $x_0\in\X$,
$Q$ is u.s.c. at $x_0$ and set $Q(x_0)$ is bounded. The
the set-valued mapping $\Phi=q+Q$ is locally bounded
around $x_0$, i.e. there exist a bounded set $W\subset\Y$
and $r>0$ such that
$$
    \Phi(x)\subseteq W,\quad\forall x\in\ball{x_0}{r}.
$$
\end{lemma}

\begin{prof}
By the continuity of $q$ at $x_0$, corresponding to $\eta>0$
there exists $r_q>0$ such that
$$
    q(x)\in\ball{q(x_0)}{\eta},\quad\forall x\in\ball{x_0}{r_q}.
$$
By the upper semicontinuity of $Q$ at $x_0$, corresponding to $\eta>0$
there exists $r_Q>0$ such that
$$
    Q(x)\subseteq\inte\ball{Q(x_0)}{\eta},\quad\forall x\in\ball{x_0}{r_Q}.
$$
Notice that, since $Q(x_0)$ is bounded, also $\ball{Q(x_0)}{\eta}$
is bounded. Thus, taking $r_\Phi=\min\{r_q,r_Q\}$, it holds
$$
   \Phi(x)=q(x)+Q(x)\subseteq\ball{q(x_0)}{\eta}+\ball{Q(x_0)}{\eta},
   \quad\forall x\in\ball{x_0}{r_\Phi},
$$
and hence it suffices to set $W=\ball{q(x_0)}{\eta}+\ball{Q(x_0)}{\eta}$.
\end{prof}

\begin{proposition}     \label{pro:exeffsol}
Suppose that the data of problem $(\mathcal{SP})$ satisfy
the following assumptions:

\noindent $(a_1)$ $(\X,\|\cdot\|)$ is of class $\UCtwo$ and $(\Y,\|\cdot\|)$
is reflexive;

\noindent $(a_2)$ $q\in\Coneone(\inte\ball{x_0}{r})$, for some
$x_0\in\X$ and $r>0$, such that $\ball{x_0}{r}\subseteq\Omega$;

\noindent $(a_3)$ $Q$ is locally closed and convex multifunction
near $x_0$ and $Q(x_0)$;

\noindent $(a_4)$ set $Q(x_0)$ is bounded and the set-valued
mapping $Q$ is u.s.c. at $x_0$;

\noindent $(a_5)$ $Q$ is metrically regular at $x_0$, for $Q(x_0)$,
with regularity modulus such that
\begin{eqnarray*}      
    \reg{Q}{x_0}< \|\der q(x_0)\|_\Lin^{-1};
\end{eqnarray*}
   
\noindent $(a_6)$  there exists $\tau>0$ such that $\Phi(\ball{x_0}{t})$
is closed for every $t\in [0,\tau]$.

\noindent Then, there exists $\epsilon_0>0$, such that for every
$\epsilon\in (0,\epsilon_0]$ problem $(\mathcal{SP}_{x_0,\epsilon})$
admits a $C$-efficient pair $(x_\epsilon,y_\epsilon)\in\fr
\ball{x_0}{\epsilon}\times\Phi(x_\epsilon)$.
\end{proposition}

\begin{prof}
Under the above hypotheses it is possible to apply Theorem
\ref{thm:PCPset}. According to it, there exists a positive $\epsilon_0$
such that for every $\epsilon\in (0,\epsilon_0]$ the image
$\Phi(\ball{x_0}{\epsilon})$ is a convex subset of $\Y$. Now, fix
any $\epsilon\in (0,\epsilon_0]$ and consider the corresponding
localized problem $(\mathcal{SP}_{x_0,\epsilon})$. Observe that
$\Phi(\ball{x_0}{\epsilon})$ is compact with respect to the weak topology
in $\Y$. Indeed, as it is norm closed and convex, it is also weakly
closed. Besides, since $Q(x_0)$ is bounded, by virtue of Lemma
\ref{lem:cuscblocb} the mapping $\Phi$ turns out to be locally bounded around $x_0$.
Thus, up to a reduction in the value of $\epsilon_0$, one can assume
that $\Phi(\ball{x_0}{\epsilon})$ is bounded. So the reflexivity of $(\Y,\|\cdot\|)$
entails that $\Phi(\ball{x_0}{\epsilon})$ is weakly compact.
By virtue of Theorem 6.5 (a) in \cite{Jahn04}, there exists an element
$y_\epsilon\in\Phi(\ball{x_0}{\epsilon})$, which is $C$-minimal.
This means that there is $x_\epsilon\in\ball{x_0}{\epsilon}$, with
$y_\epsilon\in\Phi(x_\epsilon)$, such that $(x_\epsilon,y_\epsilon)$
is a $C$-efficient pair for $(\mathcal{SP}_{x_0,\epsilon})$. 
Observe that, as a $C$-minimal element of $\Phi(\ball{x_0}{\epsilon})$,
$y_\epsilon$ must be in $\fr\Phi(\ball{x_0}{\epsilon})$.
As noted in Remark \ref{rem:mainthm} $(r_3)$, since $x_\epsilon\in
\Phi^{-1}(y_\epsilon)$, it is $x_\epsilon\in\fr\ball{x_0}{\epsilon}$.
This completes the proof.
\end{prof}

\begin{remark}
As a comment to Proposition \ref{pro:exeffsol}, it should be noted
that its thesis is trivial if $\X$ and $\Y$ are finite-dimensional
Euclidean spaces, because $\Phi(\ball{x_0}{\epsilon})$ is compact.
In an abstract space setting, under the hypotheses of Proposition
\ref{pro:exeffsol}, it is possible to state that $\ball{x_0}{\epsilon}$
is weakly compact (recall Remark \ref{rem:UCspace}($r_3$)).
Nevertheless, since $q$ may not be continuous with respect
to the weak topologies, already the set $q(\ball{x_0}{\epsilon})$ may happen
to be not weakly compact, in the absence of convexity assumptions.
\end{remark}

The next result, which comes as a further consequence of
Theorem \ref{thm:PCPset},
is an optimality condition useful for detecting solution pairs of
$(\mathcal{SP}_{x_0,\epsilon})$. It can be regarded as a scalarization
method, relying on the use of the following Lagrangian function
$L:\X\times\Y^*\longrightarrow\R\cup\{-\infty\}$
$$
    L(x,y^*)=\langle y^*,q(x)\rangle+\inf_{y\in Q(x)}\langle y^*,y\rangle,
$$
where $\Y^*$ denotes the dual space of $\Y$, whose null vector
is marked by $\nullv^*$, and $\langle\cdot,\cdot\rangle:\Y^*\times\Y
\longrightarrow\R$ denotes the canonical duality pairing $\Y^*$ with $\Y$.
To formulate such a result, one needs to consider elements
in the cone
$$
   C^+=\{y^*\in\Y^*:\ \langle y^*,y\rangle\ge 0,\quad\forall y\in C\}.
$$

\begin{proposition}    \label{pro:Lagsca}
Under the hypotheses of Proposition \ref{pro:exeffsol},
corresponding with any $\epsilon\in (0,\epsilon_0]$ and
with a $C$-efficient pair $(x_\epsilon,y_\epsilon)$,
there exists $y^*_\epsilon\in C^+\backslash\{\nullv^*\}$ such that
$x_\epsilon$ solves the scalar problem
$$
    {\rm minimize}\ L(x,y^*_\epsilon)\quad\hbox{ subject to }
   x\in\ball{x_0}{\epsilon}.
$$
\end{proposition}

\begin{prof}
By Proposition \ref{pro:exeffsol}, there exists a pair $(x_\epsilon,
y_\epsilon)\in\fr\ball{x_0}{\epsilon}\times\Phi(x_\epsilon)$, which is
$C$-efficient for problem $(\mathcal{SP}_{x_0,\epsilon})$. Since
$y_\epsilon\in\Phi(x_\epsilon)=q(x_\epsilon)+Q(x_\epsilon)$, there
exists $v_\epsilon\in Q(x_\epsilon)$ such that $y_\epsilon=
q(x_\epsilon)+v_\epsilon$. Recall that the set $\Phi(\ball{x_0}
{\epsilon})$ is closed, convex and with nonempty interior
(remember Remark \ref{rem:mainthm} $(r_3)$).
Since $y_\epsilon-C$ is convex and $(y_\epsilon-C)\cap
\Phi(\ball{x_0}{\epsilon})=\{y_\epsilon\}$, the Heidelheit theorem
applies. Consequently, there exist $y^*_\epsilon\in\Y^*
\backslash\{\nullv^*\}$ and $\alpha\in\R$ such that
\begin{eqnarray}     \label{in:linsep1}
   \langle y^*_\epsilon,y\rangle\le\alpha,\quad\forall
  y\in y_\epsilon-C
\end{eqnarray}
and
\begin{eqnarray}     \label{in:linsep2}
   \langle y^*_\epsilon,y\rangle\ge\alpha,\quad\forall
  y\in \Phi(\ball{x_0}{\epsilon}).
\end{eqnarray}
From inequality $(\ref{in:linsep1})$  it follows
\begin{eqnarray}     \label{in:linsep3}
   \langle y^*_\epsilon,y_\epsilon\rangle-
   \langle y^*_\epsilon,y\rangle\le\alpha,
   \quad\forall y\in C.
\end{eqnarray}
In particular, as it is $\nullv\in C$, one has $\langle y^*_\epsilon,
y_\epsilon\rangle\le\alpha$.
On the other hand, as $y_\epsilon\in \Phi(\ball{x_0}{\epsilon})$, then
from inequality $(\ref{in:linsep2})$, it is also $\langle y^*_\epsilon,
y_\epsilon\rangle\ge\alpha$, whence it results in
$$
    \langle y^*_\epsilon,y_\epsilon\rangle=\alpha.
$$
On account of the last equality, one sees that 
inequality $(\ref{in:linsep3})$ implies $y^*_\epsilon\in C^+$.

Now, recalling that $\Phi(x)=q(x)+Q(x)$, from inequality
$(\ref{in:linsep2})$ one obtains for every $x\in\ball{x_0}{\epsilon}$
$$
   \langle y^*_\epsilon,q(x)\rangle+\langle y^*_\epsilon,y\rangle
   \ge \langle y^*_\epsilon,q(x_\epsilon)+v_\epsilon\rangle,
   \quad\forall y\in Q(x).
$$
In particular, for $x=x_\epsilon$ it holds
$$
   \langle y^*_\epsilon,q(x_\epsilon)\rangle+\langle y^*_\epsilon,y\rangle
   \ge \langle y^*_\epsilon,q(x_\epsilon)+v_\epsilon\rangle,
   \quad\forall y\in Q(x_\epsilon).
$$ 
This allows one to deduce that
$$
    \langle y^*_\epsilon,v_\epsilon\rangle=
    \min_{y\in Q(x_\epsilon)}\langle y^*_\epsilon,y\rangle.
$$
According to the definition of $L$, one finds
$$
    L(x,y^*_\epsilon)= \langle y^*_\epsilon,q(x)\rangle+\inf_{y\in Q(x)}
   \langle y^*_\epsilon,y\rangle\ge
    \langle y^*_\epsilon,q(x_\epsilon)\rangle+\min_{y\in Q(x_\epsilon)}
  \langle y^*_\epsilon,y\rangle=L(x_\epsilon,y^*_\epsilon),\quad
    \forall x\in\ball{x_0}{\epsilon}.
$$
This completes the proof.
\end{prof}

\begin{remark}
$(r_1)$ It is worth noting that, under the assumptions of
Proposition \ref{pro:Lagsca} the Lagrangian function
$L$ can be written
$$
   L(x,y^*)=\langle y^*,q(x)\rangle+\min_{y\in Q(x)}\langle y^*,y\rangle
$$
in a neighbourhood of $x_0$. Indeed, recall that each element
$y^*\in\Y^*$ is also weakly continuous. As already seen, since
$Q(x_0)$ is bounded and $Q$ is u.s.c. at $x_0$, $Q$ turns out to
be locally bounded. Therefore, in the reflexive space $(\Y,\|\cdot\|)$
each set $Q(x)$ is weakly compact, for $x$ near $x_0$, with the
consequence that $y^*$ attains its minimum on it.

$(r_2)$ A feature of Proposition \ref{pro:Lagsca} to be commented
is that it establishes a scalarization condition which is typical in
problems with convex graph objective, even though the graph of
$\Phi$ is not necessarily convex. Indeed, according to Proposition
\ref{pro:Lagsca}, a $C$-efficient pair turns out to be a solution for
a scalar problem involving the Lagrangian function $L$, what is
more than a mere stationarity condition for $L$. This happens
by virtue of the convexity principle, which enables one to exploit
the ``hidden convexity" of the problem. In this concern, notice
that, even if the function $x\mapsto\inf_{y\in Q(x)}\langle y^*,y
\rangle$ is convex under the hypotheses of Proposition \ref{pro:Lagsca},
the function $x\mapsto L(x,y^*)$ may lose this property, owing
to the additional term $\langle y^*,q(x)\rangle$.
\end{remark}

%%%%%%%%%%%%%%%%%%%%%%%%%%%%%%%%%%%%%%%%%%%%%%%%%%%%%%%%%

{\small

\vskip.5cm

}

\end{document}